%
\magnification=\magstep1   
\input amstex
\UseAMSsymbols
\input pictex
\vsize=23truecm
\NoBlackBoxes
\parindent=18pt
  
   \font\rmk=cmr8    \font\itk=cmti8  \font\ttk=cmtt8

\font\gross=cmbx10 scaled\magstep1 

\def\mod{\operatorname{mod}}

\def\add{\operatorname{add}}
\def\Ker{\operatorname{Ker}}
  
\def\arr#1#2{\arrow <1.5mm> [0.25,0.75] from #1 to #2}

\def\T#1{\qquad\text{#1}\qquad}    
   

 \def\Rahmenbio#1%
   {$$\vbox{\hrule\hbox%
                  {\vrule%
                       \hskip0.5cm%
                            \vbox{\vskip0.3cm\relax%
                               \hbox{$\displaystyle{#1}$}%
                                  \vskip0.3cm}%
                       \hskip0.5cm%
                  \vrule}%
           \hrule}$$}        

\def\Rahmen#1%
   {\centerline{\vbox{\hrule\hbox%
                  {\vrule%
                       \hskip0.5cm%
                            \vbox{\vskip0.3cm\relax%
                               \hbox{{#1}}%
                                  \vskip0.3cm}%
                       \hskip0.5cm%
                  \vrule}%
           \hrule}}}

	\bigskip\bigskip
\centerline{\gross Minimal infinite submodule-closed subcategories.}
                    	\bigskip
\centerline{Claus Michael Ringel}     
	\bigskip\bigskip
\plainfootnote{}
{\rmk 2000 \itk Mathematics Subject Classification. \rmk 
Primary 
        16D90, 
        16G60. 
Secondary:
        16G20. 
        16G70. 
}

{\narrower \rmk Abstract.  Let $\ssize \Lambda$ be an artin algebra. 
We are going to consider full subcategories of $\ssize \mod\Lambda$
closed under finite direct sums and under submodules
with infinitely many isomorphism classes of indecomposable modules.  The main result asserts
that such a subcategory contains a minimal one and we exhibit some striking properties
of these minimal subcategories. These results have to be considered as essential
finiteness conditions for such module categories. 
\par}

	\bigskip\bigskip
Let $\Lambda$ be an artin algebra, and $\mod \Lambda$ the category of 
$\Lambda$-modules of finite length. All the subcategories to be considered will be
full subcategories of $\mod\Lambda$
closed under isomorphisms, finite direct sums and direct summands, but note that we
also consider individual $\Lambda$-modules which may not be of finite length. 
Let $\Cal C$ be a subcategory of $\mod\Lambda$. We say that $\Cal C$ is {\it finite} provided
it contains only finitely many isomorphism classes of indecomposable modules, otherwise
$\Cal C$ is said to be {\it infinite.} Of course, $\Cal C$ is said to be 
{\it submodule-closed}
provided for any module $C$ in $\Cal C$ also any submodule of $C$
belongs to $\Cal C$.
	\medskip
The aim of this paper is to study infinite submodule-closed 
subcategories of $\mod\Lambda$. A subcategory $\Cal C$ of $\mod \Lambda$ will be called 
{\it minimal infinite submodule-closed,} or (in this paper) just
{\it minimal,}
provided it is infinite and submodule-closed, and no proper subcategory of $\Cal C$ is both
infinite and submodule-closed. On a first thought, it is not at all clear whether minimal
subcategories do exist: the existence 
is in sharp contrast to the usual properties of infinite structures (recall that
in set theory, a set is infinite iff it contains proper subsets of the same cardinality).
	\medskip
{\bf Theorem 1.} {\it Any infinite submodule-closed subcategory of $\mod \Lambda$ contains
a minimal subcategory.}
	\medskip
Of course, the assertion is of interest only in case $\Lambda$ is representation-infinite.
But already the special case of looking at the category $\mod\Lambda$ itself, 
with $\Lambda$
representation-infinite, should be stressed: {\it The module category of 
any representation-infinite artin algebra has minimal subcategories.}
	\bigskip
Let $M$ be a $\Lambda$-module, not necessarily  of finite length.
We write $\Cal S_M$ for the class of finite length modules cogenerated by $M$.
This is clearly a submodule-closed subcategory of $\mod \Lambda$. (Conversely, any
submodule-closed subcategory $\Cal C$ of $\mod \Lambda$ is of this form: take for $M$
the direct sum of all modules in $\Cal C$, one from each isomorphism class; 
or else, it is sufficient to take just indecomposable modules in $\Cal C$.). 
	\medskip
{\bf Theorem 2.} {\it Let $\Cal C$ be a minimal subcategory of $\mod\Lambda$.
Then
\item{\rm(a)} For any natural number $d$, there are only finitely many isomorphism classes
of modules in $\Cal C$ of length $d$.
\item{\rm(b)} Any module in $\Cal C$ is isomorphic to a submodule of an indecomposable
module in $\Cal C$.
\item{\rm(c)} There is an infinite sequence of indecomposable modules $C_i$ in $\Cal C$ 
with proper inclusions
$$
 C_1 \subset C_2 \subset \cdots \subset C_i \subset C_{i+1} \subset \cdots
$$
such that also the union $M = \bigcup_i C_i$ is indecomposable and then
$\Cal C = \Cal S_M.$}
	\bigskip
As we have mentioned, Theorem 1 asserts, in particular, that the module category of 
any representation-infinite artin algebra has a minimal subcategory $\Cal C$,
and the assertion (c) of Theorem 2 yields arbitrarily large indecomposable modules in
$\Cal C$. This shows that we are in the realm of the first Brauer-Thrall conjecture
(formulated by Brauer and Thrall around 1940 and proved by Roiter in 1968): 
any representation-infinite artin algebra has indecomposable modules of arbitrarily
large length. The proof of Roiter and its combinatorial interpretation by Gabriel
are the basis of the Gabriel-Roiter measure on $\mod\Lambda$, see [R1, R2].
Using it, we have shown in [R1] that the module category of a representation-infinite
artin algebras always has a so-called take-off part: this is an infinite 
submodule-closed subcategory with property (a) of Theorem 2, 
and there is an infinite inclusion
chain of indecomposables such that also the union $M$ is indecomposable, 
as in property (c) of Theorem 2.
However,  $\Cal S_M$
usually will be a proper subcategory of the take-off part, and then
the take-off part cannot be minimal. Of course, we can apply Theorem 1
to the take-off part in order to obtain a minimal subcategory inside the take-off part.
The important feature of
the minimal categories is the following: we deal with a countable set of 
indecomposable modules which are strongly interlaced as the assertions (b) and (c) of
Theorem 2 assert. Typical examples to have in mind are the 
infinite preprojective components of hereditary algebras (see section 4).  
	\bigskip
The proof of theorem 1 will be given in section 2, the proof of theorem 2 in section 3.
These proofs depend on the  Gabriel-Roiter measure for $\Lambda$-modules,
as discussed in [R1,R2]. 
The remaining section 4 provides examples. First, we will mention
some procedures for obtaining submodule-closed subcategories.
Then, following Kerner-Takane, we will show that the preprojective component 
of a representation-infinite connected hereditary algebra $\Lambda$ 
is always a minimal subcategory. In case $\Lambda$ is tame, this is the only one, but
for wild hereditary algebras, there will be further ones.
	\bigskip
{\bf Acknowledgment.} The results have been announced at the Annual meeting of the German
Mathematical Society, Bonn 2006 and in further lectures at various occasions. In particular,
two of the Selected-Topics lectures [R3, R4] in Bielefeld were devoted to this theme.
The author is grateful to many mathematicians for comments concerning the presentation. 
	\bigskip\bigskip
{\bf 2. Proof of Theorem 1.}
	\medskip
 Given a class $\Cal X$ of modules of finite length 
(or of isomorphism classes of modules),
we denote by $\add \Cal X$ the smallest subcategory containing $\Cal X$.
We denote by $\Bbb N = \Bbb N_1$ the natural numbers starting with 1.
	\medskip
The proof will be based on results concerning the Gabriel-Roiter measure
for $\Lambda$-modules, see [R1, R2]; 
the Gabriel-Roiter measure $\mu(M)$ of a module $M$ will be 
considered either as a finite set $I$ of 
natural numbers, or else as the rational number $\sum_{i\in I} 2^{-i},$
whatever is more suitable. 
Given a Gabriel-Roiter measure $I$, let $\Cal C(I)$ be the set of isomorphism
classes of indecomposable objects in $\Cal C$ with Gabriel-Roiter measure $I$. 
An obvious adaption of one of the main results of [R1] asserts:
	\medskip
{\it There is an infinite sequence of  Gabriel-Roiter measures $I_1 < I_2 < \cdots$
such that $\Cal C(I_t)$ is non-empty for any $t\in \Bbb N$ and such that 
for any $J$ with $\Cal C(J) \neq \emptyset$, either $J = I_t$ for some $t$ or else
$J > I_t$ for all $t$. Moreover, all the sets $\Cal C(I_t)$ are finite.} (Note that
the sequence of measures $I_t$ depends on $\Cal C$, thus one should write
$I_t^{\Cal C} = I_t$; the papers [R1,R2] were dealing only with 
the case $\Cal C = \mod\Lambda$, but the proofs carry over to the more general
case of dealing with a submodule-closed subcategory $\Cal C$).
	\medskip
Since $\add \bigcup_{t\in \Bbb N} \Cal C(I_t)$ is an infinite submodule-closed 
subcategory of $\Cal C$, we may assume that
$\Cal C = \add \bigcup_{t\in \Bbb N} \Cal C(I_t)$.
In order to construct a minimal subcategory $\Cal C'$, we will construct a sequence of subcategories
$$
 \Cal C = \Cal C_0 \supseteq \Cal C_1 \supseteq \Cal C_2 \supseteq \cdots
$$
with the following properties: \newline
(a) Any subcategory $\Cal C_i$ is infinite and submodule-closed,
$$
 \Cal C_i(I_t) = \Cal C_t(I_t)    \T{for} t \le i. \tag{b}
$$
(c) If $\Cal D \subseteq \Cal C_i$ is infinite and submodule-closed, then
$$
 \Cal D(I_t) = \Cal C_t(I_t) \T{for} t \le i.  
$$

We start with $\Cal C_0 = \Cal C$ (the $t$ in conditions (b) and (c) satisfies $t \ge 1$,
thus nothing has to be verified). Assume, we have constructed $\Cal C_{i}$ for some $i\ge 0$,
satisfying the conditions (a), and the conditions (b), (c) for all pairs $(i,t)$ with
$t \le i.$ We are going to construct $\Cal C_{i+1}.$

Call a subset $\Cal X$ of $\Cal C_{i}(I_{i+1})$ {\it good}, provided there is
a subcategory $\Cal D_{\Cal X}$ of $\Cal C_{i}$ which is infinite and submodule-closed
and such that $\Cal D_{\Cal X}(I_{i+1}) = \Cal X.$ For example $\Cal C_{i}(I_{i+1})$ itself
is good (with $\Cal D_{\Cal X} = \Cal C_{i}).$ Since $\Cal C_{i}(I_{i+1})$ is a finite set,
we can choose a minimal good subset $\Cal X' \subseteq \Cal X$. For $\Cal X'$, there is
an infinite and submodule-closed subcategory $\Cal D_{\Cal X'}$ of $\Cal C_{i}$
such that $\Cal D_{\Cal X'}(I_{i+1}) = \Cal X'.$ (Note that in general neither $\Cal X'$ 
nor $\Cal D_{\Cal X'}$ will be uniquely determined: usually, there may be several possible
choices.) Let $\Cal C_{i+1} = \Cal D_{\Cal X'}$. By assumption, $\Cal C_{i+1}$ is
infinite and submodule-closed, thus (a) is satisfied. In order to show (b) for all
pairs $(i+1,t)$ with $t \le i+1$, we first consider some $t\le i$. We can
apply (c) for $\Cal D = \Cal C_{i+1} \subseteq \Cal C_i$ and see that $\Cal D(I_t) =
\Cal C_t(I_t)$, as required. But for $t = i+1,$ nothing has to be shown.
Finally, let us show (c). Thus let $\Cal D \subseteq \Cal C_{i+1}$ be an infinite
submodule-closed subcategory. Since $\Cal D \subseteq \Cal C_i$, we know by induction that
$\Cal D(I_t) = \Cal C_t(I_t)$ for $t \le i.$ It remains to show that
$\Cal D(I_{i+1}) = \Cal C_{i+1}(I_{i+1}).$ Since $\Cal D \subseteq  \Cal C_{i+1},$
we have $\Cal D(I_{i+1}) \subseteq \Cal C_{i+1}(I_{i+1}).$ But if this would be
a proper inclusion, then $\Cal X = \Cal D(I_{i+1})$ would be a good subset of
$\Cal C_{i}(I_{i+1})$ which is properly contained in $\Cal C_{i+1}(I_{i+1}) 
= \Cal D_{\Cal X'}(I_{i+1}),$ a contradiction to the minimality of $\Cal X'.$ 
This completes the inductive construction of the various $\Cal C_i.$
	\medskip
Now let
$$
 \Cal C' =  \bigcap\nolimits_{i\in \Bbb N} \Cal C_i.
$$
Of course, $\Cal C'$ is submodule-closed. Also, we see immediately
$$
 \Cal C'(I_t) = \Cal C_t(I_t)    \T{for all} t, \tag{b$'$}
$$
since $\Cal C'(I_t) = \bigcap\nolimits_{i\ge t} \Cal C_i(I_t) = \Cal C_i(I_t)$, according to 
(b).

First, we show that $\Cal C'$ is
infinite. Of course, $\Cal C'(I_1) \neq \emptyset,$
since $I_1 = \{1\}$ and a good subset of $\Cal C_0(I_1)$ has to contain at least one
simple module.  Assume that $\Cal C'(I_s) \neq \emptyset$
for some $s$, we want to see that there is $t > s$ with $\Cal C'(I_{t}) \neq \emptyset$.
For every Gabriel-Roiter measure $I$, let $n(I)$ be the minimal number $n$ with $I \subseteq 
[1,n]$, thus $n(I)$ is the length of the modules in $\Cal C(I)$. Let $n(s)$ be the maximum
of $n(I_j)$ with $j \le s$, thus $n(s)$ is the maximal length of the modules in 
$\bigcup_{j\le s} \Cal C(I_j).$ Let $s'$ be a natural number such that
$n(I_j) > n(s)pq$ for all $j > s'$ (such a number exists, since the modules in $I_j$ with
$j$ large, have large length); here $p$ is the maximal length of an indecomposable projective
module, $q$ that of an indecomposable injective module.

We claim that $\Cal C'(I_j) \neq \emptyset$ for some $j$ with $s < j \le s'.$ Assume for the
contrary that $\Cal C'(I_j) = \emptyset$ for all $s < j \le s'.$ We consider $\Cal C_{s'}.$
Since $\Cal C_{s'}$ is infinite, there is some $t > s$ with $\Cal C_{s'}(I_t) \neq \emptyset,$
and we choose $t$ minimal. 
Now for $s < j \le s'$, we know that $\Cal C_{s'}(I_j) = \Cal C_{j}(I_j) = 
\Cal C'(I_j) = \emptyset,$ according to (b) and (b$'$). This shows that $t > s'.$ 
Let $Y$ be an indecomposable module with isomorphism class in $\Cal C_{s'}(I_t)$. 
Let $X$ be a Gabriel-Roiter submodule of $Y$. Then $X$ belongs to $\Cal C_{s'}(I_j)$
with $j < t$. If $j \le s,$ then the length of $X$ is bounded by $n(s)$, and therefore
$Y$ is bounded by $n(s)pq$ (see [R2], 3.1 Corollary), in contrast to the fact that 
$n(I_t) > n(s)pq.$ This is the required contradiction. Thus $\Cal C'$ is infinite.

Now, let $\Cal D$ be an infinite submodule-closed subcategory of $\Cal C'$.
We show that $\Cal D[I_t] = \Cal C'[I_t]$ for all $t$. Consider some fixed $t$ and choose an 
$i$ with $i \ge t$.
Since $\Cal C' \subseteq \Cal C_i$, we see that $\Cal D[t] = \Cal C_t[t]$ the given $t$,
according to (b) for $\Cal C_i$. But according to (b$'$), we also
know that $\Cal C'[t] = \Cal C_t[t]$. This completes the proof. 
	\bigskip\bigskip

{\bf 3. Proof of Theorem 2.}
	\medskip
	
We refer to [R1] for the proof of (a) and for the 
construction of an inclusion chain 
$$
 C_1 \subset C_2 \subset \cdots \subset C_i \subset C_{i+1} \subset \cdots
$$
with indecomposable union, as asserted in (c). In [R1] these assertions have been
shown for the take-off part of $\mod\Lambda$, but the same proof with only minor
modifications, carries over to minimal categories. 
	\medskip
To complete the proof of (c), we only have to note the following:
By construction, $\Cal S_M$ contains all the modules $C_i$, thus $\Cal S_M$ is not
finite. But of course, $\Cal S_M \subseteq \Cal C$. Namely, if $X$ is a finite length module
which is cogenerated by $M$, then there are finitely many maps $f_i\:X \to M$ such that
the intersection of the kernels is zero. But there is some $j$ such that the images of 
all the maps $f_i$ are contained in $C_j$, therefore $X$ is cogenerated by $C_j$ and
thus belongs to $\Cal C$. The minimality of $\Cal C$ implies that $\Cal S_M = \Cal C.$
	\medskip
It remains to proof part (b) of Theorem 2.
We will need some general observations which may be of
independent interest.  Recall that a module is
said to be of 
{\it finite type,} provided it is the direct sum of (may-be infinitely many) copies of 
a finite number of modules of finite length). 
	\medskip

(1) {\it If $\Cal S_M$ is minimal, then $M$ is not of finite type.} 
	\medskip
Proof: Assume that $M$ is of finite type, let $M_1,\dots,M_t$ be the indecomposable
direct summands of $M$, one from each isomorphism class. We may
assume that they are indexed with increasing Gabriel-Roiter measures, thus
$\mu(M_i) \le \mu(M_j)$ for $i \le j.$ 
Let $M'$ be the direct sum of all indecomposable modules in $\Cal S_M$ which are
not isomorphic to $M_t$. Since $\Cal S_M$ is infinite, also
$\Cal S_{M'}$ is infinite, and of course $\Cal S_{M'} \subseteq \Cal S_M.$
Assume that $M_t$ belongs to $\Cal S_{M'}$. Then $M_t$ is cogenerated by a finite
number of indecomposable modules $N_1,\dots,N_s$ which are direct summands of $M'$.
Thus, we have an embedding $u: M_t \to N$, where $N$ is a direct sum of
copies of these modules $N_i$. 
Also, the modules $N_i$ are cogenerated by $M$, thus there is
an embedding $u'\:N \to M^r$ for some $r$. Altogether, $u'u\: M_t \to M^r$. 
Now $\mu(M_t) = \max_{1\le i \le t} \mu(M_i)$, and therefore $u'u$ is a split
monomorphism. Consequently, also $u\:M_t \to N$ is a split monomorphism, and therefore
one of the modules $N_i$ is isomorphic to $M_t$, in contrast to our construction.
	\medskip
(2) {\it If $\Cal S_M$ is minimal and $M' \subseteq M$ is a cofinite submodule, then
$\Cal S_{M'} = \Cal S_M.$}
	\medskip
Proof: Of course, $\Cal S_{M'} \subseteq \Cal S_M.$ Since we assume that $\Cal S_M$
is minimal, we only have to show that $\Cal S_{M'}$ is infinite. Assume, for the contrary,
that $\Cal S_{M'}$ is finite. This implies that $M'$ is of finite type (see [R5]),
say
$M' = \bigoplus_{i\in I} M'_i,$  so that the modules $M'_i$ belong to only finitely many
isomorphism classes.
Let $U$ be a submodule of $M$ of finite length such that $M'+U = M.$ Now $M'\cap U$ is
a submodule of $M'$ of finite length, thus it is contained in some
$M' = \bigoplus_{i\in J} M'_i,$ where $J$ is a finite subset of $I$. It follows that
$M = U \oplus  \bigoplus_{i\in I\setminus J} M'_i,$ and this again is a module of
finite type. But this contradicts (1). 
	\medskip
(3) {\it Assume that $\Cal C = \Cal S_M$ is minimal and let $M_0$ be a submodule of $M$
of finite length. If $X$ belongs to $\Cal C$, then there is an embedding $u\:X \to M$
such that $M_0\cap u(X)  = 0.$}
	\medskip
Proof. Let $X$ be of finite length and cogenerated by $M$.
We want to construct inductively maps $f\:X \to M$ such that $M_0 \cap f(X)= 0$ and
such that the length of $\Ker(f)$ decreases. 
As start, we take as $f$ the zero map. The process will end when $\Ker(f) = 0.$

Thus, assume that we have given some map $f\:X \to M$ with $M_0\cap f(X) = 0$ and
$\Ker(f) \neq 0.$  We are going to construct
a map $g\:X \to M$ such that first $M_0\cap g(X) = 0$ and second,
$\Ker(g)$ is a proper submodule of $\Ker(f)$. 
Let $M_1 = M_0+f(X),$  this is a submodule of finite length of $M$.
Choose a submodule $M'$ of $M$ with $M_1\cap M' = 0,$ and maximal with this property.
Note that $M'$ is a cofinite submodule of $M$ (namely,
$M/M'$ embeds into the injective hull of $M_1$, and with $M_1$ also its
injective hull has finite length). 
According to (2), we know that $\Cal S_{M'} = \Cal S_M
= \Cal C$, thus $X$ belongs to $\Cal S_{M'}$. This means that $X$ is cogenerated by $M'$.
In particular, since $\Ker(f) \neq 0$, there is a map $f'\:X \to M'$ such that
$\Ker(f)$ is not contained in $\Ker(f')$. Let $g = (f,f')\:X \to M_1\oplus M' \subseteq M$.
Then $\Ker(g) = \Ker(f)\cap \Ker(f')$ is a proper submodule of $\Ker(f).$ Also,
the image $g(X)$ is contained in $f(X)+f'(X) \subseteq f(X)+M'$. Since $M_1 + M' =
M_0\oplus f(X)\oplus M'$, we see that $M_0\cap g(X) = 0.$ 

This completes the induction step. After finitely many steps, 
we obtain in this way an embedding $u$ of $X$ into $M$ such that $u(X)\cap M_0 = 0.$ 
	\medskip
(3$'$) {\it Assume that $\Cal C = \Cal S_M$ is minimal. If $X, Y $ are submodules of $M$ of finite
length, then also $X\oplus Y$ is isomorphic to a submodule of $M$.}
	\medskip
Proof: If $X, Y $ are submodules of $M$, then $X\oplus Y$ is cogenerated by $M$.
	\medskip
(3$''$) {\it Assume that $\Cal C = \Cal S_M$ is minimal. If $C$ belongs to $\Cal C$, 
then the direct sum of countably many copies of $C$ can be embedded into $M$.}
	\medskip
Proof: Assume, there is given an embedding $u_t\: C^t\to M$, where $t \ge 0$ is a natural number.
Let $M_0 = u_t(C^t).$ According to (3), we find an embedding $u\:C \to M$ such that
$M_0\cap u(C) = 0.$ Thus, let $u_{t+1} = u_t\oplus u\:C^{t+1} = C^t\oplus C \to M.$
	\bigskip
Proof of part (b) of Theorem 2.
Let $C$ be a module in $\Cal C$. Let $M = \bigcup_i C_i$ be 
as constructed in (c), thus all the $C_i$ are indecomposable and $\Cal S_M = \Cal C$.
According to (3), there is an embedding $u\:C \to M$. Now the image of $u$ lies in some
$C_i$, thus $u$ embeds $C$ into the indecomposable module $C_i$.
	\medskip
At least one consequence of Theorem 2 (b) should be mentioned. 
If $S$ is a simple $\Lambda$-module, write  $[X\:S]$ for the Jordan-H\"older
multiplicity of $S$ in the $\Lambda$-module $X$. 
	\medskip
{\bf Corollary.} {\it Let $\Cal C$ be a minimal subcategory. For any natural number
$d$, there is an indecomposable module $C$ in $\Cal C$ with the following property:
if $S$ is a simple $\Lambda$-module with $[Y\:S]\neq 0$ for some $Y$
in $\Cal C$, then $[C\:S] \ge d$.}
	\medskip
Proof: We consider the simple $\Lambda$-modules $S$ 
such that there exists a module $Y(S)$
in $\Cal C$ with $[Y(S)\:S] \neq 0,$ and let $Y = \bigoplus Y(S)$
where the summation extends over all isomorphism classes of such simple modules $S$.
Given a natural numer $d$, let us consider $Y^d$. According to assertion (b) of
Theorem 2, there is an indecomposable $\Lambda$-module $C$ such that $Y^d$
embeds into $C$. But this implies that $[C\:S] \ge [Y^d\:S] =d[Y\:S] \ge d[Y(S):S] 
\ge d.$
	\medskip
Note that the corollary provides a strengthening of the assertion of the
first Brauer-Thrall conjecture: {\it A representation-infinite artin algebra
has indecomposable representations $X$ such that all non-zero Jordan-H\"older
multiplicities of $X$ are arbitrarily large.}

	\bigskip\bigskip
{\bf 4. Examples.} 
	\medskip
Let us start to mention some ways for obtaining submodule-closed subcategories $\Cal C$
of $\mod\Lambda$.
	\medskip
\item{$\bullet$} Of course, we can consider the module category $\mod\Lambda$ itself. 
\item{$\bullet$} If $\Cal I$
is a two-sided ideal of $\Lambda$, then the $\Lambda$-modules annihilated by 
$\Cal I$ form a submodule-closed subcategory
(this subcategory is just the category of all $\Lambda/\Cal I$-modules).
\item{$\bullet$} As we have mentioned in section 3, we may start with an 
arbitrary (not necessarily finitely generated) module $M$, and consider the
subcategory $\Cal S_M$ of all finite length modules cogenerated by $M$. This
subcategory $\Cal S_M$ is submodule closed, and any submodule-closed subcategory of 
$\mod \Lambda$ is obtained in this way. 
\item{$\bullet$} The special case of dealing with
$M = {}_\Lambda\Lambda$ has been studied often in representation theory; the
modules in $\Cal S_{{}_\Lambda\Lambda}$ are called the {\it torsionless}
$\Lambda$-modules. Artin algebras with $\Cal S_{{}_\Lambda\Lambda}$ finite 
have 
quite specific properties, for example their representation dimension is bounded by 3.
\item{$\bullet$} 
 The categories $\Cal A(<\! \gamma)$ and $\Cal A(\le\! \gamma)$ of all modules $X$ in
$\Cal A = \mod \Lambda$ with
Gabriel-Roiter measure $\mu(X) < \gamma$, or  $\mu(X) \le \gamma$, respectively;
here $\gamma \in \Bbb R$ and $\mu$ is the Gabriel-Roiter measure 
(or a weighted Gabriel-Roiter measure). 
\item{$\bullet$} In particular, the take-off subcategory of $\mod\Lambda$ (as 
introduced in [R1]) is submodule-closed (and it is infinite iff $\Lambda$ 
is representation-infinite). 
\item{$\bullet$} If $\Lambda$ has global dimension $n$, then the subcategory $\Cal C$ of all
modules of projective dimension at most $n-1$ is closed under cogeneration (and extensions)
(this is mentioned for example in [HRS], Lemma II.1.2.).
	\medskip
Given such a submodule-closed subcategory $\Cal C$, one may ask whether it is finite or
not, and in case it is infinite, it should be of interest 
to look at the corresponding minimal subcategories.
	\bigskip

{\bf Example 1 (Kerner-Takane).} {\it Let $\Lambda$ be a connected hereditary artin algebra
of infinite representation type. The preprojective component of $\mod\Lambda$
is a minimal subcategory.}
	\medskip
Proof. Kerner-Takane ([KT], Lemma 6.3.)
have shown: For every $b\in \Bbb N$, there is $n = n(b)\in \Bbb N$ with the following
property: If $P, P'$ are indecomposable projective modules, then $\tau^{-i}P'$ is cogenerated by 
$\tau^{-j}P$, for all $0 \le i \le b$ and $n \le j$. Assume that
$\Cal C$ is the additive subcategory
given by an infinite set of indecomposable preprojective modules. We claim that the cogeneration
closure of $\Cal C$ contains all the preprojective modules $X$. Indeed, let
$X = \tau^{-b}P'$ with $P'$ indecomposable projective. 
Choose a corresponding $n(b)$. Since $\Cal C$ contains infinitely many
isomorphism classes of indecomposable preprojective modules, there is some $C = \tau^{-j}P$ in $\Cal C$
with $n \le j$ and $P$ indecomposable projective. According to Kerner-Takane, $X$ is cogenerated by
$C$.
	\bigskip
{\bf Example 2.} {\it Any tame concealed algebras $\Lambda$ has a unique minimal 
subcategory $\Cal C$, namely the subcategory of all preprojective
modules.}
	\medskip
We need the following two well-known results:
	\smallskip
{\bf Lemma 1.} {\it 
Let $P$ be preprojective with defect $\delta(P) = -1$, let $R$ be indecomposable regular
with regular radical $R'$ (this means that $R/R'$ is simple regular). Assume that $f\:P \to R$
is a map with image not contained in $R'$. Then $f$ is a monomorphism or an epimorphism.}
	\smallskip
Proof: Assume that $f$ is not a monomorphism, let $K$ be its kernel and $I$ its image.
Then $K$ is preprojective,
in particular $\delta(K) \le -1,$ and therefore $\delta(P/K) \ge 0.$ But $I$ as a 
submodule of $R$ with $\delta(I) \ge 0$ has to be regular. Since $I$ is not included in $R'$,
it follows that $I = R.$
	\smallskip
{\bf Lemma 2.} {\it 
Let $P$ be preprojective with defect $-d$. Then there are $d$ preprojective modules $P_i$
of defect $-1$, and surjective maps $f_i\:P \to P_i$ such that $(f_i)\:P \to \bigoplus P_i$ is
injective.}
	\smallskip
Proof: 
There is an exact sequence $0 \to P @>f>> G^d @>>> Y \to 0$, where $G$ is the generic module,
and $Y$ is a direct sum of Pr\"ufer modules (see [RR]). Let $P_i$ be the image of the
composition of $f$ and the $i$-the canonical projection $G^d \to G$, and let
$f_i$ be the corresponding map $P \to G$ with image $P_i$. 
	\smallskip

Now, let $\Cal C$ be an infinite submodule-closed subcategory of $\mod \Lambda$,
where $\Lambda$ is a tame concealed algebra. We want to show that $\Cal C$ contains all
the preprojective modules. Recall that for $M$ an indecomposable 
$\Lambda$-module, $-6 \le \delta(M) \le 6.$

If $\Cal C$ contains infinitely many indecomposable preinjective modules, it also contains
arbitrarily large preinjective modules with defect 1 (by the dual of Lemma 2): any indecomposable
preinjective module $Q$ has a preinjective submodule $Q'$ of defect $1$ such that
$|Q'| \ge \frac16 |Q|.$ The dual of Lemma 1 asserts that a preinjective module $Q$ of defect -1
has a regular submodule $R$ with $|R| \ge |Q|-e$, where $e$ is the maximum of the length of the
simple regular modules in an exceptional tube (or $e = 1$ in case $\Lambda$ has only two simple modules).
Next, assume that $\Cal C$ contains infinitely many indecomposable regular modules. If they are
of bounded length, then Brauer-Thrall 1 yields arbitrarily large indecomposable modules $M$
cogenerated by these regular modules, and these modules $M$ have to be preprojective. If $\Cal C$
contains large indecomposable regular modules, then also large preprojective modules, by Lemma 2.
Altogether, we see that $\Cal C$ contains infinitely many preprojective modules. But for every
natural number $n$ there is $n'$ such that any indecomposable preprojective module of length at least $n'$
will cogenerate all the preprojective modules of length at most $n$. This shows that $\Cal C$ contains
all the preprojective modules. --- On the other hand, the subcategory of all preprojective modules
is infinite and closed under cogeneration. This completes the proof.
	\medskip
{\bf Remark.} Preprojective components are always submodule closed, but an infinite
preprojective component $\Cal P$ does not have to be minimal. First of all, 
$\Cal P$ may contain indecomposable injective modules, whereas this cannot happen for
a minimal subcategory, as the part (b) of Theorem 2 shows. But also preprojective
components without indecomposable injective modules may not be minimal. For example,
consider the algebra with quiver
$$
{\beginpicture
\setcoordinatesystem units <1cm,1cm>
\multiput{$\circ$} at 0 0  2 0  4 0 /
\arr{1.7 0}{0.3 0}
\arr{3.7 0.1}{2.3 0.1}
\arr{3.7 -.1}{2.3 -.1}
\put{$a$} at 0 .3
\put{$b$} at 2 .3
\put{$c$} at 4 .3
\setdots <1mm>
\setquadratic
\plot 1 0.2  2 0.6  3 0.3 /
\endpicture}
$$
with one zero relation. Then the preprojective component $\Cal P$ contains 
indecomposables which are faithful, but also countable many indecomposables
$X$ with $X_a = 0.$ Clearly, the subcategory of $\Cal P'$ of all modules
$P$ in $\Cal P$ with $P_a = 0$ is a proper subcategory which is both
infinite and submodule closed (and actually, $\Cal P'$ is minimal). 

	\bigskip

{\bf Example 3.} Let $\Cal I$ be a twosided ideal in $\Lambda$. The category of $\Lambda$-modules
annihilated by $\Cal I$ is obviously submodule-closed and of course equivalent (or even equal)
to the category of all $\Lambda/\Cal I$-modules. If $\Lambda/\Cal I$ is representation-infinite, then
$\mod \Lambda/\Cal I$ will contain a minimal subcategory. Consider for
example the generalized Kronecker-algebra $K(3)$ with three arrows $\alpha,\beta,\gamma.$
The one-dimensional ideals of $K(3)$ correspond bijectively to the elements of the
projective plane $\Bbb P^2$, say $a = (a_0:a_1:a_2) \in \Bbb P^2$ yields the ideal
$\Cal I_a = \langle a_0\alpha+a_1\beta+a_2\gamma \rangle$. Let $\Cal C_a$ be 
the additive subcategory
of $\mod K(3)$ of all preprojective
$K(3)/\Cal I_a$-modules. Then these are pairwise different 
subcategories (the intersection of any two of these subcategories is the subcategory of 
semisimple projective modules). In particular, {\it if the base field is infinite, there are
infinitely many subcategories in $\mod K(3)$ 
which are minimal.} (Note that
the preprojective $K(3)$-modules provide a further minimal subcategory.)
	\medskip
The minimal subcategories exhibit here can be distinguished by looking at the
corresponding annihilators (the annihilator of a subcategory $\Cal C$ is the ideal of
all the elements $\lambda\in \Lambda$ which annihilate all the modules in $\Cal C$).
The next example will show that usually there are also different minimal
subcategories which have the same annihilator. Note that a submodule-closed
subcategory $\Cal C$ has zero annihilator if and only if all the projective
modules belong to $\Cal C.$ 
	\bigskip
{\bf Example 4.} Here is an artin algebra $\Lambda$ with 
different minimal categories containing all the
indecomposable projective modules. Consider the hereditary algebra $\Lambda$
with quiver $Q$
$$
{\beginpicture
\setcoordinatesystem units <1cm,1cm>
\multiput{$\circ$} at 0 0  2 0  4 0 /
\arr{1.7 0.1}{0.3 0.1}
\arr{1.7 -.1}{0.3 -.1}
\arr{3.7 0.1}{2.3 0.1}
\arr{3.7 -.1}{2.3 -.1}
\put{$a$} at 0 .3
\put{$b$} at 2 .3
\put{$c$} at 4 .3
\put{$\alpha$} at 1 0.3
\put{$\alpha'$} at 1 -.3
\put{$\beta$} at 3  .3
\put{$\beta'$} at 3 -.3
\endpicture}
$$
We denote by $Q_{ab}$ the full subquiver $Q_{ab}$
with vertices $a,b$, by $Q_{bc}$ that with vertices $b,c$.

As we know, the preprojective component $\Cal C$ of $\mod\Lambda$
is a minimal subcategory. Of
course, it contains all the projective $\Lambda$-modules, but 
it contains also, for example, the indecomposable $\Lambda$-module $X$ with dimension
vector $(3,2,0);$ note that the restriction of $X$ to $Q_{ab}$ 
is indecomposable and neither projective nor semisimple. 

Second, let $\Cal D$ be the full subcategory of $\mod \Lambda$ consisting of
all the $\Lambda$-modules such that the restriction to $Q_{ab}$
is projective 
and the restriction to $Q_{bc}$ is preprojective. 
Clearly, $\Cal D$ is submodule-closed, and it
is obviously infinite: If $Y$ is a $\Lambda$-module with $Y_a = 0,$ 
define $\underline Y$ as follows:
the restrictions of $Y$ and $\underline Y$ 
to $Q_{bc}$ should coincide, whereas the restriction of
$\underline Y$ to $Q_{ab}$ should be a direct sum
of indecomposable projectives of length 3; in particular, $\underline Y_a =
Y_b^2.$ By $Y \mapsto \underline Y$ we obtain an embedding of the category of  
preprojective 
Kronecker modules into $\Cal D,$ which yields all the indecomposable modules
in $\Cal D$ but the simple projective one. It follows easily that $\Cal D$
is minimal. Of course, $\Cal D \neq \Cal C$, and note that also $\Cal D$ 
contains all the projective $\Lambda$-modules. 

We can exhibit even a third minimal subcategory which contains all
the projective $\Lambda$-modules, by looking at the
full subcategory $\Cal E$ of $\Lambda$-modules such that the restriction to
$Q_{ab}$ is the direct sum of a projective and a semisimple module, whereas
the restriction to $Q_{bc}$ is projective. Again, clearly 
$\Cal E$ is submodule-closed.
In order to construct an infinite family of indecomposable modules in $\Cal E$,
we use covering theory: The following quiver is part of the 
universal cover $\widehat Q$ of $Q$  
$$
{\beginpicture
\setcoordinatesystem units <.7cm,.8cm>
\multiput{$1$} at  0 0  2 0  3 0  5 0  6 0  8 0  9 0  11 0 
  2.5 2  5.5 2  8.5 2    1 1
  10 1 /
\multiput{$2$} at    4 1  7 1   /
\arr{0.7 0.7}{0.3 0.3}
\arr{1.3 0.7}{1.7 0.3}
\arr{3.7 0.7}{3.3 0.3}
\arr{4.3 0.7}{4.7 0.3}
\arr{6.7 0.7}{6.3 0.3}
\arr{7.3 0.7}{7.7 0.3}
\arr{9.7 0.7}{9.3 0.3}
\arr{10.3 0.7}{10.7 0.3}

\arr{2.1 1.7}{1.4 1.3}
\arr{2.9 1.7}{3.6 1.3}
\arr{5.1 1.7}{4.4 1.3}
\arr{5.9 1.7}{6.6 1.3}
\arr{8.1 1.7}{7.4 1.3}
\arr{8.9 1.7}{9.6 1.3}
\multiput{$\beta$} at 1.6 1.7  4.6 1.7  7.6 1.7 /
\multiput{$\beta'$} at 3.5 1.7  6.5 1.7  9.5 1.7 / 
\multiput{$\alpha$} at 0.2 0.6  3.2 0.6  6.2 0.6  9.2 0.6  / 
\multiput{$\alpha'$} at 1.9 0.65  4.9 0.65  7.9 0.65  10.9 0.65 /

\endpicture}
$$
and the numbers inserted form the dimension vector of a lot of indecomposable
modules $M$. If we require in addition 
that the maps $\alpha$ and $\alpha'$ starting at
the same vertex have equal kernels, then there is a unique
isomorphism class $M = Y_3$ with this dimension vector.
In a similar way, we can construct for any
natural number $n$ an indecomposable representation $Y_n$ of $\widehat Q$ of
length $2+5n$ (with top of length $n$). The kernel condition 
assures that the $\Lambda$-module
which is covered by $M = Y_3$, or more generally, by $Y_n$,
belongs to $\Cal E$ (note that the kernel condition
means that the restriction of $M$ to any subquiver of type $\widetilde {\Bbb D}_4$
has socle of length 3).
If $\Cal E'$ is a minimal subcategory inside $\Cal E$, then $\Cal E'$
is different from $\Cal C$ and $\Cal D.$
(Remark: The $\Lambda$-module covered by $Y_1$ is indecomposable projective and has
Gabriel-Roiter measure $(1,3,7)$, this is the measure  $I_3$ for $\lambda$. 
One may show that the $\Lambda$-module covered by $Y_2$ has
Gabriel-Roiter measure $(1,3,7,12)$ and that this is the measure
$I_4$. For $t\ge 5$, the measures $I_t$ are not yet known; it would be
interesting to decide whether the intersection of the take-off part of $\mod\Lambda$
and $\Cal E$ is infinite or not.)

	\bigskip\bigskip
{\bf References.}
	\medskip
\item{[HRS]} D.Happel, I.Reiten, S.Smal{\o}: Tilting in Abelian Categories and Quasitilted Algebras.
   Memoirs AMS 575 (1996)
\item{[KT]} O.Kerner, M.Takane: Mono orbits, epi orbits and elementary vertices of 
 representation infinite quivers. 
 Comm. Alg. 25, 51-77 (1997),
\item{[RR]} I.Reiten, C.M.Ringel: Infinite dimensional representations of canonical algebras. 
   Canadian Journal of Mathematics 58 (2006), 180-224. 
\item{[R1]} C.M.Ringel: The Gabriel-Roiter measure.
 Bull. Sci. math. 129 (2005), 726-748.
\item{[R2]} C.M.Ringel:  Foundation of the Representation Theory of Artin Algebras, Using the Gabriel-Roiter Measure. 
 In: Trends in Representation Theory of Algebras and Related Topics (ed: 
 de la Pena and Bautista). Contemporary Math. 406. Amer.Math.Soc. (2006), 105-135. 
\item{[R3]}  C.M.Ringel: Minimal infinite cogeneration-closed subcategories.
Selected Topics, Bielefeld 2006. www.mathematik.uni-bielefeld.de/$\sim$sek/select/minimal.pdf
\item{[R4]}  C.M.Ringel: Pillars.
Selected Topics, Bielefeld 2006. \newline
 www.mathematik.uni-bielefeld.de/$\sim$sek/select/pillar.pdf
\item{[R5]} C.M.Ringel: The first Brauer-Thrall conjecture. 
In: Models, Modules and Abelian Groups. In Memory of A. L. S. Corner. Walter de Gruyter, Berlin 
 (ed: B. Goldsmith, R. G\"obel) (2008), 369-374. 
\item{[R6]}  C.M.Ringel: Gabriel-Roiter inclusions and Auslander-Reiten theory. J.Algebra (to appear)

	\bigskip\bigskip

{\rmk Fakult\"at f\"ur Mathematik, Universit\"at Bielefeld \par
POBox 100\,131, \ D-33\,501 Bielefeld, Germany \par
e-mail: \ttk ringel\@math.uni-bielefeld.de \par}

	\bigskip\bigskip

\bye